\documentclass{birkjour}
\usepackage{amsthm,amsfonts,amssymb,amsmath,enumerate,cite}

 \newtheorem{thm}{Theorem} 

  \newtheorem{prop}[thm]{Proposition}
  \theoremstyle{definition}
  
\newtheorem{ex}[thm]{Example}
  \theoremstyle{remark}
  
\DeclareMathOperator{\BV}{BV}
\DeclareMathOperator{\CE}{CE}
\DeclareMathOperator{\Hom}{Hom}
\DeclareMathOperator{\Sh}{Sh}

\def\Ai{A_\infty}
\def\BVi{\BV_\infty}
\def\Li{L_\infty}
\def\ri{r_\infty}
\def\S{\widehat{S}}
\def\B{\mathcal{B}}
\def\Z{\mathbb Z}
\def\g{\mathfrak g}

\newcommand{\abs}[1]{{\lvert#1\rvert}}

\begin{document}

 %
 %
 %
 %
 %
 %
 %
 %
 \title[$\ri$-Matrices, Triangular $\Li$-Bialgebras and
 Quantum$_\infty$ Groups]{$\ri$-Matrices, Triangular $\Li$-Bialgebras\\
   and Quantum$_\infty$ Groups}

\author[D. Bashkirov]{Denis Bashkirov}
\address {School of Mathematics\\University of Minnesota\\
  Minneapolis, MN 55455\\USA}
\email{bashk003@umn.edu}

\thanks{This work was supported by the World Premier International
  Research Center Initiative (WPI Initiative), MEXT, Japan, the
  Institute for Mathematics and its Applications with funds provided
  by the National Science Foundation, and a grant from the Simons
  Foundation (\#282349 to A.~V.).}

\author[A. A. Voronov]{Alexander A. Voronov}
\address {School of Mathematics\\University of Minnesota\\
  Minneapolis, MN 55455\\USA\\
and\\
Kavli IPMU (WPI)\\ University of
  Tokyo\\Kashiwa, Chiba 277-8583\\Japan}
\email{voronov@umn.edu}



 \date{February 1, 2015}

 \begin{abstract}
   A homotopy analogue of the notion of a triangular Lie bialgebra is
   proposed with a goal of extending basic notions of the theory of
   quantum groups to the context of homotopy algebras and, in
   particular, introducing a homotopical generalization of the notion
   of a quantum group, or quantum$_\infty$-group.
 \end{abstract}

 \maketitle

\section{Introduction}

\subsection{Conventions and Notation}

We will work over a ground field $k$ of characteristic zero. A
differential graded (dg) vector space $V$ will mean a complex of
$k$-vector spaces with a differential of degree one. The degree of a
homogeneous element $v \in V$ will be denoted by $\abs{v}$. In the
context of graded algebra, we will be using the Koszul rule of signs
when talking about the graded version of notions involving symmetry,
including commutators, brackets, symmetric algebras, derivations,
\emph{etc}., often omitting the modifier \emph{graded}. For any
integer $n$, we define a \emph{translation} (or \emph{$n$-fold
  desuspension$)$} $V[n]$ of $V$: $V[n]^p := V^{n+p}$ for each $p \in
\Z$. For two graded vector spaces $V$ and $W$, we define grading on
the space $\Hom(V,W)$ of $k$-linear maps $V \to W$ by $\abs{f} := n-m$
for $f \in \Hom(V^m,W^n)$.

\subsection{Quantum Groups}

Recall that a \emph{quantum group} in the sense of Drinfeld and Jimbo
is an associative, coassociative Hopf algebra $A$ subject to the
condition of being \textit{quasitriangular} \cite{Dr1}. The latter
implies, in particular, the existence of a solution $\mathcal{R}$ to
the \emph{quantum Yang-Baxter equation}
$\mathcal{R}^{12}\mathcal{R}^{13}\mathcal{R}^{23}=\mathcal{R}^{23}\mathcal{R}^{
  13}\mathcal{R}^{12}$ set up in $A$. More conceptually, the
quasitriangularity condition provides data needed to put a braided
structure on the monoidal category of left $A$-modules.

The most basic examples of quantum groups appear as
\emph{quantizations} or certain types of \emph{deformations} (in the
sense of Hopf algebras) of universal enveloping algebras and algebras
of functions on groups. In the first case, starting with a Lie algebra
$\g$ and a Hopf-algebra deformation $U_h(\g)$ of its universal
enveloping algebra $U(\g)$, one passes to the
``\emph{$($semi$)$classical limit}''
$\delta(x):=\frac{\Delta_h(x)-\Delta_h^{op}(x)}{h}$ thus equipping
$U(\g)$ with a \emph{co-Poisson-Hopf structure} with $\delta: U(\g)
\to U(\g) \otimes U(\g)$ being the \emph{co-Poisson cobracket}. In
particular, the restriction $\delta|_{\g}$ becomes a well-defined
cobracket on $\g$ turning it into a Lie bialgebra.

\subsection{Quantization of Triangular Lie Bialgebras}
\label{classical}

The above process can be reversed: as it was shown in \cite{Re}, any
fi\-nite-di\-men\-sion\-al Lie bialgebra $(\g,[,],\delta)$ can be
\textit{quantized}, meaning that one can always come up with a
Hopf-algebra deformation $U_h(\g)$ whose classical limit, in the sense
of the above formula, agrees with $\delta$. While a priori $U_h(\g)$
is just a Hopf algebra, one would really be interested in having a
quasitriangular structure on it.  As a special case, it was shown in
\cite{Dr83} that such a structure exists, when $\g$ is a
\textit{triangular} Lie bialgebra. This class of Lie bialgebras is
defined as follows: let $\g$ be a Lie algebra and $r \in \g\otimes\g$
(``\emph{a classical $r$-matrix}'') be a skew-symmetric element
satisfying the \emph{classical Yang-Baxter equation}
\[
[r_{12},r_{13}]+[r_{12},r_{23}]+[r_{23},r_{13}] = 0,
\]
which can be conveniently restated in the form of the
\emph{Maurer-Cartan equation}
\[
[r,r] = 0
\]
taking place in the graded Lie algebra $S(\g[-1])[1]$ with respect to
the \emph{Schouten bracket} for elements $r$ of degree one: $r \in
(S(\g[-1])[1])^1 = (S(\g[-1]))^2 = S^2(\g[-1])[2] = \g \wedge
\g$. Such an element $r$, called a \emph{Maurer-Cartan element}, gives
rise to a Lie cobracket on $\g$ in the form of the coboundary
$\partial_{\CE}(r): \g \to \g \wedge \g$ of $r$ taken in the cochain
\emph{Chevalley-Eilenberg complex} of $\g$ with coefficients in $\g
\wedge \g$ (here, $r$ is regarded as a 0-cocycle). The compatibility
with the Lie-algebra structure on $\g$ is packed into the relation
$\partial_{\CE}^2(r)=0$, thus guaranteeing that $\g$ with such a
cobracket is indeed a Lie bialgebra. The \emph{co-Jacobi identity},
which could be rewritten as
\begin{equation}
\label{co-Jacobi}
[\partial_{\CE} (r), \partial_{\CE} (r)] = 0,
\end{equation}
follows from the following statement, which is an odd version of the
Hamiltonian correspondence in Poisson geometry, if one regards
$S(\g[-1])[1]$ as the shifted Gerstenhaber algebra of functions on the
odd Poisson manifold $(\g[-1])^*$ and $\Hom (\g[-1], S(\g[-1]))$ as
the graded Lie algebra of vector fields on $(\g[-1])^*$.
\begin{prop}
  \label{Hamiltonian}
  The linear map
\[
\partial_{\CE}: S(\g[-1])[1] \to \Hom (\g[-1], S(\g[-1]))
\]
is a graded Lie-algebra morphism.
\end{prop}

\noindent
A \textit{triangular Lie bialgebra} is a Lie algebra $\g$ provided
with a Lie cobracket $\partial_{\CE}(r)$ coming out of an $r$-matrix
$r$.  A basic statement concerning this class of Lie bialgebras is
that $U(\g)$ can be quantized to a \textit{triangular Hopf algebra}
$U_h(\g)$. This condition is stronger than being quasitriangular, and
in particular, the category of (left) modules over a triangular Hopf
algebra turns out to be symmetric monoidal, as opposed to just being
braided.

\subsection{The Homotopy Quantization Program}

The upshot of the above construction is that there is a source of
quantum groups coming from the data of a Lie algebra $\g$ and a
solution of the Maurer-Cartan equation in $S(\g[-1])[1]$. The goal of
our project is to promote this construction to the realm of homotopy
Lie algebras.  In particular, this would generalize the work
\cite{BSZ} done for the case of Lie 2-bialgebras. Here is an outline
of our program:

\begin{enumerate}[1.]
\item Develop the notion of a \textit{triangular $\Li$-bialgebra}
  extending the classical one. In analogy with the classical case, the
  input data for this construction consists of an $\Li$-algebra $\g$
  and a solution $r$ of a generalized Maurer-Cartan equation set up in
  an appropriate algebraic context;
\item Show that the universal enveloping algebra $U(\g)$ of a
  triangular $\Li$-bialgebra $\g$ admits a natural \textit{homotopy
    co-Poisson-Hopf} structure;
\item Extend the \emph{Drinfeld twist} construction \cite{Dr1}, which
  equips a cocommutative Hopf algebra with a new, triangular
  coproduct, to the homotopical context. Apply it to the case of
  universal enveloping algebra $U(\g)$ of the previous step to obtain
  a \emph{quantum$_\infty$ group}.
\end{enumerate}

\noindent
The current paper is dedicated to describing the first step of the
construction, which we believe might be interesting on its own.

The second step is work in progress. While the universal enveloping
algebra $U(\g)$ of an $\Li$-algebra $\g$ is a strongly homotopy
associative (or $\Ai$-) algebra that also turns out to be a
cocommutative, coassociative coalgebra object in Lada-Markl's, see
\cite{lada-markl}, symmetric monoidal category of $\Ai$-algebras, we
would be interested in verifying that $U(\g)$ is actually a
\emph{Hopf$_\infty$ algebra}, that is, possesses an antipodal map
satisfying certain compatibility conditions. We would also need to
translate an $\Li$-bialgebra structure on an $\Li$-algebra $\g$ into a
co-Poisson$_\infty$ structure on the Hopf$_\infty$ algebra
$U(\g)$. This would result in providing $U(\g)$ with the structure of
a cocommutative co-Poisson$_\infty$ coalgebra object in Lada-Markl's
symmetric monoidal category of $\Ai$-algebras.

Furthermore, we would like to develop deformation theory of homotopy
Hopf algebras and use it to quantize homotopy Lie bialgebras. A
different approach to quantization of homotopy Lie bialgebras (using
the framework of PROPs) was taken in \cite{merk:qu}, in which a
different notion of a homotopy Hopf algebra (or rather, homotopy
bialgebra) was used. That notion depends on the choice of a minimal
resolution of the bialgebra properad. The notion we outline above
appears to be more canonical.

In the future we would also be interested in investigating what this
program produces for $\Li$-algebras arising in the geometric context,
such as generalized Poisson geometry, $\Li$-algebroids, or
$\BVi$-geometry.

\section{The Big Bracket and $\Li$-Bialgebras}
\label{BB}

Recall that the structure of a (strongly) homotopy Lie algebra (or
\textit{$\Li$-algebra}) on a graded vector space $\g$ may be given by
a \emph{codifferential}, \emph{i.e}., a degree-one, square-zero
coderivation $D$ such that $D(1) = 0$, on the graded cocommutative
coalgebra $S(\g[1])$ equipped with the shuffle comultiplication.  The
date given by $D$ is equivalent to a collection of ``higher Lie
brackets'' $l_k: S^k(\g[1]) \to \g[1]$, $k\geq 1$, of degree one
obtained by restriction of $D$ to the $k$th symmetric component of
$S(\g[1])$ followed by projection to the cogenerators. The condition
$D^2 = 0$ is equivalent to the \emph{higher Jacobi identities},
homotopy versions of the Jacobi identity. Outside of deformation
theory, nontrivial examples of homotopy Lie algebras are known to
arise in the context of multisymplectic geometry \cite{rogers,
  baez-hoffnung}, Courant algebroids \cite{roytenberg}, and closed
string field theory \cite{zwiebach,ksv,markl-loop}.

In order to discuss the structure of an $\Li$-bialgebra on a graded
vector space $\g$, we need to mix the graded Lie algebra $\Hom
(S(\g[1]),\g[1])$, used to define the structure of an $\Li$-algebra on
$\g$, with the graded Lie algebra $\Hom (\g[-1]), S (\g[-1]))$, used
to define the cobracket on $\g$ in the classical, Lie algebra setting
in Section~\ref{classical}. Consider the graded vector space
\[
\B : = \prod\limits_{m,n\geq 0} \Hom (S^m(\g[1]),S^n(\g[-1]))[2]
\]
and provide it with the structure of a graded Lie algebra given by the
graded commutator, called the \emph{big bracket},
\[
[f,g] : = f \circ g - (-1)^{\abs{f}\cdot \abs{g}} g \circ f
\]
under the \emph{circle}, or $\cup_1$ \emph{product}, cf.\
\cite{gerstenhaber-voronov} and \cite{terilla:qdt}:
{\small \begin{multline*}
  (f \circ g) (x_1 \dots x_n) \\ := \sum_{\sigma \in \Sh_{k,l}}
  (-1)^\varepsilon f(x_{\sigma(1)} \dots x_{\sigma(k)}
  g(x_{\sigma(k+1)} \dots x_{\sigma(n)})_{(1)}) \, g(x_{\sigma(k+1)}
  \dots x_{\sigma(n)})_{(2)},
\end{multline*}}
where $x_1, \dots, x_n \in \g[1]$,
\begin{gather*}
f  \in \prod\limits_{m\geq 0} \Hom (S^{k+1}(\g[1]),S^m(\g[-1]))[2],\\
g \in \prod\limits_{m\geq 0} \Hom (S^{l}(\g[1]),S^m(\g[-1]))[2],
\end{gather*}
$n = k+l$ --- otherwise we set $(f \circ g) (x_1 \dots x_n) = 0$,
$\Sh_{k,l}$ is the set of $(k,l)$ \emph{shuffles}: permutations
$\sigma$ of $\{1,2, \dots, n\}$ such that $\sigma(1) < \sigma(2) <
\dots < \sigma (k)$ and $\sigma(k+1) < \sigma(k+1) < \dots < \sigma
(n)$, $\varepsilon = \abs{x_\sigma} + \abs{g}(\abs{x_{\sigma(1)}} +
\dots + \abs{x_{\sigma(k)}})$, $(-1)^{\abs{x_\sigma}}$ is the
\emph{Koszul sign} of the permutation of $x_1 \dots x_n$ to
$x_{\sigma(1)} \dots x_{\sigma(n)}$ in $S(\g[1])$, and we use
Sweedler's notation to denote the result $g_{(1)} \otimes g_{(2)}$ of
applying to $g \in S(\g[-1])[2]$ the (shifted) comultiplication
$S(\g[-1])[2] \to S(\g[-1])[2] \otimes S(\g[-1])$ followed by the
projection $S(\g[-1])[2] \to \g[1]$ onto the cogenerators in the first
tensor factor. This graded Lie algebra $\B$, under the assumption that
$\dim \g < \infty$ and in a slightly different incarnation, was
introduced by Y.~Kosmann-Schwarzbach \cite{kos-sch1} in relation to
Lie bialgebras and later used by O.~Kravchenko \cite{kravchenko:sh-bv}
in relation to $\Li$-bialgebras. The graded Lie algebra $\B$ has the
property that its Maurer-Cartan elements represent $\Li$ brackets and
cobrackets on $\g$, as well as mixed operations, comprising the
structure of an $\Li$-bialgebra on $\g$. Here we adopt Kravchenko's
approach and define an \emph{$\Li$-bialgebra structure} on $\g$ as a
Maurer-Cartan element $\mu$ in the subalgebra
\[
\B^+ := \prod\limits_{m,n\geq 1} \Hom (S^m(\g[1]), S^n(\g[-1]))[2]
\]
of the graded Lie algebra $\B$.\footnote{Maurer-Cartan elements in
  $\B$ would correspond to more general, \emph{curved}
  $\Li$-bialgebras.} This means
\begin{gather}
\nonumber
\mu = \sum\limits_{m,n\geq 1} \mu_{mn},\\
\label{mu_mn}
\mu_{mn}: S^m(\g[1]) \to S^n(\g[-1])[2] \quad \text{of degree 1},\\
\nonumber
[\mu, \mu] = 0.
\end{gather}

\section{$\ri$-Matrices and Triangular $\Li$-Bialgebras}

For an $\Li$-algebra $\g$, one can generalize the Schouten bracket to
an $\Li$ structure on $S(\g[-1])[1]$ by extending the higher brackets
$l_k$ on $\g$ as graded multiderivations of the graded commutative
algebra $S(\g[-1])$. This $\Li$ structure may also be described via
\emph{higher derived brackets} (in the semiclassical limit) on the
$\BVi$-algebra $S(\g[-1])$, see \cite[Example 3.4]{BaVo}. The $\Li$
structure can be naturally extended to the completion
\[
\S(\g[-1])[1] := \prod\limits_{n\geq 0} S^n(\g[-1])[1].
\]

While investigating the deformation-theoretic meaning of solutions $r
= r(\lambda) \in \lambda \S(\g[-1])[1][[\lambda]]$, where $\lambda$ is
the \emph{deformation parameter}, that is to say, a (degree-zero)
formal variable, of the \emph{generalized Maurer-Cartan equation}
\begin{align}
\label{MCE}
l_1(r) + \frac{1}{2!} l_2(r \odot r) + \frac{1}{3!} l_3(r \odot r
\odot r) + \dots = 0,
\end{align}
where $\odot$ refers to multiplication in $S(V)$ for $V =
\S(\g[-1])[2]$, certain analogies can be drawn with basic constructions
of the theory of quantum groups.

Note that an $\Li$-algebra structure on $\g$ endows the graded Lie
algebras $\B$ and $\B^+$ with the structure of a dg Lie
algebra. Namely, bracketing with the Maurer-Cartan element $l_1 + l_2
+ \dots \in \prod\limits_{m \geq 1} \Hom (S^m(\g[1]), \g[-1])[2]$,
representing the $\Li$-algebra structure on $\g$, creates a
differential
\[
d \gamma : = [ l_1 + l_2 + \dots, \gamma]
\]
on $\B$ and $\B^+$, compatible with the ``big'' bracket.

An $\Li$-\emph{morphism} $\varphi$ from the $\Li$-algebra
$\S(\g[-1])[1]$ to the dg Lie algebra $\B^+$, that is to say, a
morphism $S(\S(\g[-1])[2])\to S(\B^+[1])$ of dg coalgebras mapping 1
to 1, amounts to defining a series of degree-zero linear maps
$\varphi_n: S^n(\S(\g[-1])[2]) \to \B^+[1]$ for all $n \ge 1$
satisfying the following compatibility conditions:
{\small\begin{multline*}
  d \varphi_n (x_1 \odot \dots \odot x_n) \\
 + \frac{1}{2} \sum_{k= 1}^{n-1} \sum_{\sigma \in \Sh_{k,n-k}}
  (-1)^\varepsilon [\varphi_k(x_{\sigma(1)} \odot \dots \odot
  x_{\sigma(k)}), \varphi_{n-k} (x_{\sigma(k+1)} \odot \dots \odot
  x_{\sigma(n)})]\\ =\sum_{m=1}^n \sum_{\tau \in \Sh_{m,n-m}}\!\!
  (-1)^{\abs{x_\tau}} \varphi_{n-m+1}( l_m(x_{\tau(1)} \odot \dots
  \odot x_{\tau(m)}) \odot x_{\tau(m+1)} \odot \dots \odot
  x_{\tau(n)}),
\end{multline*}}
where $x_1, \dots , x_n \in \S(\g[-1])[2]$ and $\varepsilon =
{\abs{x_\sigma} + \abs{x_{\sigma(1)}} + \dots +
  \abs{x_{\sigma(k)}}}$. There is a \emph{canonical $\Li$-morphism}
$\varphi: \S(\g[-1])[1] \to \B^+$, which may be defined by the maps
\begin{align*}
  \varphi_n(x_1 \odot \dots \odot x_n)(y) := l_{n + p} (x_1 \odot
  \dots \odot x_n \odot N(y)),
\end{align*}
where $x_1, \dots , x_n \in \S(\g[-1])[2]$, $y \in S^{p}(\g[1])$,
$N(y) \in S^{p}(\S(\g[-1])[2])$, $p \ge 1$, and $N: S(\g[1]) \to
S(\S(\g[-1])[2])$ is the graded-algebra morphism induced by the
obvious linear map $\g[1] \hookrightarrow \S(\g[-1])[2]
\hookrightarrow S(\S(\g[-1])[2])$. The following theorem generalizes
Proposition~\ref{Hamiltonian} to the $\Li$ setting.

\begin{thm}
  The above maps $\varphi_n$, $ n \ge 1$, define an $\Li$-morphism
\[
\varphi:  \S(\g[-1])[1] \to \B^+
\]
from the $\Li$-algebra $\S(\g[-1])[1]$ to the dg Lie algebra $\B^+$.
\end{thm}

\noindent
The proof of the theorem is a straightforward checkup that reduces the
statement to the higher Jacobi identities for the $\Li$ brackets on
$\S(\g[-1])[1]$. This theorem also generalizes Kravchenko's result
\cite[Theorem~19]{kravchenko:sh-bv}, which provides an $\Li$-morphism
from an $\Li$-algebra $\g$ to the graded Lie algebra $\Hom (\g, \g)$.

An \textit{$\ri$-matrix} $r$ is a \emph{$($generalized$)$
  Maurer-Cartan element} $r= r(\lambda)$ in the $\Li$-algebra $\lambda
\S(\g[-1])[1][[\lambda]]$, \emph{i.e}., a degree-one solution of the
generalized Maurer-Cartan equation \eqref{MCE}. Sending an
$\ri$-matrix to the subalgebra $\B^+$ of the big-bracket dg Lie
algebra $\B$ under an $\Li$-morphism $\varphi: \S(\g[-1])[1] \to \B^+$
would yield a Maurer-Cartan element $\mu' = \mu'(\lambda)$:
\[
\mu' := \varphi (e^r) = \varphi_1(r) + \frac{1}{2!}
\varphi_2(r \odot r) + \frac{1}{3!} \varphi_3 (r \odot r \odot r) +
\dots,
\]
depending on the deformation parameter $\lambda$, in the dg Lie
algebra $\lambda\B^+[[\lambda]]$:
\[
d \mu' + \frac{1}{2} [\mu', \mu'] = 0,
\]
or, equivalently, a Maurer-Cartan element
\[
\mu = \mu' + l_1 + l_2 + \dots
\]
in the graded Lie algebra $\lambda\B^+[[\lambda]]$:
\[
[\mu, \mu] = 0,
\]
giving rise to an $\Li$-bialgebra structure on $\g$, as per
Section~\ref{BB}, in analogy with $\partial_{\CE}(r)$ giving rise to
an ordinary Lie cobracket in the classical, nonhomotopical case, see
also Example~\ref{example} below.

We call the $\Li$-bialgebra $(\g, \mu)$ produced out of an
$\Li$-algebra $\g$ and an $\ri$-matrix $r$ by transferring it to a
Maurer-Cartan element $\mu' = \varphi (e^r)$ in $\B^+$ via the
canonical $\Li$-morphism $\varphi$, as above, a \emph{triangular
  $\Li$-bialgebra}.

\begin{ex}
\label{example}
In the case of a classical Lie algebra $\g$, the graded Lie algebra
$\S(\g[-1])[1]$ is just the completed graded Lie algebra of
(right-)invariant multivector fields on the corresponding local Lie
group with the Schouten bracket (or, equivalently, up to degree shift,
functions on the formal odd Poisson manifold $(\g[-1])^*$), and the
canonical $\Li$-morphism $\varphi$ is just linear: $\varphi =
\varphi_1 = \partial_{\CE}$, equal to the Chevalley-Eilenberg
differential of 0-cochains of the Lie algebra $\g$ with coefficients
in the graded $\g$-module $\S(\g[-1])[2]$. The fact that $\varphi$ is
an $\Li$-morphism translates into being a dg Lie morphism,
\emph{i.e}., satisfying two compatibility conditions
\begin{align*}
  d\varphi (x) &= 0,\\
  [\varphi(x), \varphi(y)] &= \varphi([x,y]),
\end{align*}
where $d$ is the operator taking the big bracket with the Lie
structure $l_2 \in \Hom (S^2(\g[1]), \g[-1])[2] \subset \B^+$:
\[
d \alpha := [l_2, \alpha].
\]
The first condition means that $\partial_{\CE} (\varphi (x))
= \partial_{\CE}^2(x) = 0$, and the second states that
$\partial_{\CE}$ is a Lie-algebra morphism, which is the assertion of
Proposition~\ref{Hamiltonian}. An $r_\infty$-matrix is a solution $r
\in (\S(\g[-1])[1])^1 = (S^2(\g[-1]))^2 = \g \wedge \g$ of the
generalized Maurer-Cartan equation \eqref{MCE}, which turns into the
classical Maurer-Cartan equation $[r,r] = 0$ in this case, and thereby
$r$ is just a classical $r$-matrix. Thus, the transfer $\varphi(e^r) =
\varphi (r)$ of an $r$-matrix is a Lie cobracket $\varphi (r)$ on
$\g$, satisfying the compatibility condition with the Lie bracket and
the co-Jacobi identity \eqref{co-Jacobi}, resulting in the structure
of a triangular Lie bialgebra on $\g$. Here we ignored the deformation
parameter $\lambda$, because $\S(\g[-1])[1]$ is just a graded Lie
algebra and the generalized Maurer-Cartan equation in $\S(\g[-1])[1]$
and the morphism $\varphi$ have only finitely many terms.
\end{ex}

\begin{ex}
  When $\g$ is a dg Lie algebra, the picture is dramatically different
  from the classical picture of Example~\ref{example}. Now an
  $r_\infty$-matrix $r$ may have many more components than just one in
  $\g \wedge g$:
\[
r \in (\S(\g[-1])[1])^1 = (\S(\g[-1]))^2 = \prod\limits_{n\geq 0}
(S^n(\g[-1]))^2 = \prod\limits_{n\geq 1} (S^n(\g[-1]))^2.
\]
However, $\S(\g[-1])[1]$ is just a dg Lie algebra, and the generalized
Maurer-Cartan equation \eqref{MCE} is still classical,
\[
dr + \frac{1}{2} [r,r] = 0.
\]
The $\Li$-morphism $\varphi$ is still linear $\varphi = \varphi_1$,
\emph{i.e}., $\varphi$ is a dg Lie-algebra morphism. The transfer
$\mu' = \varphi(e^r) = \varphi(r)$ of $r$ via $\varphi$ will result in
a $\Li$-bialgebra structure on $\g$ with higher operations $\mu_{mn}$,
see \eqref{mu_mn}, which are trivial for all pairs $(m,n)$ but those
with $m=1$, $n \ge 1$ and $m=2$, $n = 1$. Thus, even in the case of a
dg Lie algebra $\g$, our construction creates a triangular
$\Li$-bialgebra, rather than just a (dg) triangular Lie algebra.
\end{ex}


\subsection*{Acknowledgment}
The authors are grateful to Yvette Kosmann-Schwarzbach for useful
remarks. A.~V. also thanks IHES, where part of this work was done, for
its hospitality.

\bibliographystyle{plain}

\end{document}